\documentclass[12pt,a4paper,oneside]{article}
\usepackage[T1]{fontenc}
\usepackage{graphicx}						%Einbinden von Grafiken
\usepackage{amsmath, amsthm, amssymb}
\usepackage{stmaryrd}						%bessere Formeln
\usepackage{amstext}						%Text in Formeln verwenden
\usepackage{amsfonts}
\usepackage[mathscr]{eucal}
\usepackage{enumerate}
\usepackage[shortlabels]{enumitem}
\usepackage{tikz}
\usetikzlibrary{arrows,shapes}
\usetikzlibrary{mindmap,trees}
\usetikzlibrary{shapes.misc}
\usetikzlibrary{decorations.pathmorphing, patterns}
\usepackage{pgfplots}
\pgfplotsset{compat=newest}

\makeatletter

\pgfdeclarepatternformonly[\LineSpace]{my north east lines}{\pgfpointorigin-0.1pt}{\pgfqpoint{\LineSpace}{\LineSpace}}{\pgfqpoint{\LineSpace}{\LineSpace}}%
{
	\pgfsetcolor{\tikz@pattern@color}
	\pgfsetlinewidth{0.4pt}
	\pgfpathmoveto{\pgfqpoint{0pt}{\LineSpace}}
	\pgfpathlineto{\pgfqpoint{100pt}{\LineSpace}}
	\pgfusepath{stroke}
}
\makeatother
\makeatletter
\let\@fnsymbol\@arabic
\makeatother
\newdimen\LineSpace
\tikzset{
	line space/.code={\LineSpace=#1},
	line space=8pt
}

\newcommand{\A}{A_{p,q}^s(\mathbb{R}^{n})}
\newcommand{\Asp}{A_{p,q}^s}

\newcommand{\As}[1]{A_{p,q}^{#1}(\mathbb{R}^{n})}
\newcommand{\AsP}[1]{A_{p,q}^{#1}}

\newcommand{\Cinf}{C^{\infty}(\mathbb{R}^n \times (0,T))}

\newcommand{\Fou}{\mathcal{F}}

\newcommand{\K}[1]{G_{#1}}
\newcommand{\La}[3]{L_{#1}((0,T),\,#2,\,#3)}

\newcommand{\Lloc}{L_1^{\text{loc}}(\Rn)}

\newcommand{\Lp}[1]{L_{#1}(\mathbb{R}^{n})}
\newcommand{\Lv}{L_v((0,T),\,b,\,X)}

\newcommand{\N}{\mathbb{N}}

\newcommand{\R}{\mathbb{R}}
\newcommand{\Rn}{\mathbb{R}^{n}}

\newcommand{\Rx}{\mathbb{R}^{n+1}}
\newcommand{\Rxt}{\mathbb{R}^{n+1}_{+}}

\newcommand{\Sm}[1]{W^\alpha_{#1}}
\newcommand{\Sw}{S(\mathbb{R}^{n})}
\newcommand{\Sd}{S'(\mathbb{R}^{n})}

\newcommand{\T}{T_{u_0}}

\newcommand{\bal}{\begin{align}}
\newcommand{\eal}{\end{align}}
\newcommand{\als}{\begin{align*}}
\newcommand{\las}{\end{align*}}
\newcommand{\eq}{\begin{eqnarray}}
\newcommand{\eqs}{\begin{eqnarray*}}
\newcommand{\hra}{\hookrightarrow}
\newcommand{\qe}{\end{eqnarray}}
\newcommand{\qes}{\end{eqnarray*}}

\newcommand{\dt}{\partial_t}

\newcommand{\bli}{\begin{list}{}{\labelwidth1.7em\leftmargin2.1em}}
\newcommand{\eli}{\end{list}}

\theoremstyle{definition}

\newtheorem{rem}{Remark}[section]

\theoremstyle{plain}

\newtheorem{prop}{Proposition}[section]
\newtheorem{theorem}{Theorem}[section]

\setlist[enumerate,1]{label={\upshape(\roman*)}}

\providecommand{\keywords}[1]
{
	\small	
	\textbf{\textit{Keywords and Phrases.}} #1
}
\providecommand{\classification}[1]
{
	\small	
	\textbf{\textit{Math Subject Classifications.}} #1
}
\begin{document}
\title{\Large{On the Cauchy problem for a semi-linear hyperdissipative heat equation}}
\date{\today}
\author{Franka Baaske\thanks{\emph{Corresponding author.} University of Applied Sciences Mittweida, Faculty Applied Computer Sciences \& Biosciences, Technilumplatz 17, 09648 Mittweida, Germany. Email: baaske@hs-mittweida.de}\hspace{.3cm} and \hspace{.1cm} Romaric Kana Nguedia\thanks{Friedrich Schiller University Jena, Faculty of Mathematics and Computer Sciences, Ernst-Abbe Platz 2, 07743 Jena, Germany. Email:  romaric.kana.nguedia@uni-jena.de}}
\maketitle
\begin{center}
	\textbf{Abstract}
\end{center}	
The paper is concerned with the Cauchy problem for a semi-linear hyperdissipative heat equation in Besov and Triebel-Lizorkin spaces which is related to the generalized Gauss-Weierstrass semi-group via Duhamel's principle.\\ Using caloric smoothing properties of the semi-group we prove existence and uniqueness of mild and strong solutions which are local in time. Moreover, we study well-posedness of the problem. \\[1ex]
\noindent
\classification{46E35, 35K25, 35K55, 35Q35}\\[1ex]

\tableofcontents%\clearpage
%

%}
%
%%\researchsupported{Please insert information concerning research grant support here.}
%%% If there is an additional footnote on page 1, place ``\makethankshere'' subsequent
%%% to that footnote and use the class option ``nothanks''.
%
%\acknowledgments{We would like to thank Hans Triebel and Hans-Jürgen Schmeißer for useful comments and suggestions.}

%\dedicated{Dedicated to Professor Vakhtang Kokilashvili}

\section{Introduction}
We consider the Cauchy problem for a semi-linear hyperdissipative heat equation
\begin{align}\label{ghsup}
\begin{array}{rcll}
\dt u(x,t) + (-\Delta_x)^\alpha u(x,t)&=&|u(x,t)|^{r-1} u(x,t), & x\in\Rn,\;0<t<T,\\
u(x,0)&=&u_0(x),  & x\in\Rn 
\end{array}      
\end{align}
where $r\geq2$, $0<T\leq \infty$, $\alpha\in\N$ and $2\leq n\in\N$. We use standard notation. In particular $\Delta_x$ denotes the Laplacian with respect to the space variable and $\partial_t=\partial/\partial t$ stands for the time derivative.The interest in higher natural powers of the Laplacian and this particular semi-linear inhomogeneity originates in the mathematical description of self-organization of cells and cellular pattern formation in developmental biology, see e.g. \cite{BGM}, \cite{DHBG}, \cite{DHDWBG} and \cite{Tri17}. In general, one wishes to model the amount of a kind of material depending on the flux of the material and its net rate of production. The former is expressed in terms of the Laplacian where its higher powers are used to describe turbulent instabilities observed  e.g. in dense bacterial suspensions. The semi-linear term models the structure or pattern formation of the material, for instance polar alignment if $r=3$. Concerning general powers $(-\Delta)_x^{\alpha}$, $\alpha>0$ we refer in the context of generalized heat and Navier-Stokes equations for instance to \cite{BaS22},  \cite{BST}, \cite{KS},\cite{M}, \cite{MG} and \cite{MYZ}.  Due to the biological background we ask for real valued solutions (under the assumption that initial values $u_0$ are real).\\
Our aim is to prove existence and uniqueness of mild and strong solutions in weighted Lebesgue spaces with respect to the Bochner integral $L_v((0,T),\,b,\,X)$ defined in \eqref{Lv} and \eqref{Lvinf} below.
We follow the method developed in \cite{Tri13} and consider a solution of \eqref{ghsup} as fixed point of the operator $T_{u_0}$ given as
\begin{align}\label{T0}
T_{u_0}u(x,t):=\Sm{t}u_0(x)+\int_0^t \Sm{t-\tau}|u|^{r-1}u(x,\tau)\text{d}\tau,\quad x\in\Rn,\;0<t<T.
\end{align}
Here $\Sm{t}\omega$ with $\omega\in\Sd$ is defined on the Fourier side as
\eq\label{Salph}
\left(\Sm{t}\omega\right)^{\wedge}(\xi):=e^{-t|\xi|^{2\alpha}}\omega^{\wedge}(\xi),\quad\xi\in\Rn,\;t>0,
\qe
and refers to the generalized Gauss-Weierstrass semi-group as introduced in \cite{BaS19}.% with $\alpha=2$.
This can be reformulated as
\eq
\Sm{t}\omega(x)=(\K{t}\ast\omega)(x),\quad x\in\Rn
\qe
where
\eq\label{K}
\K{t}(x)=(2\pi)^{-n/2}\left(e^{-t|\xi|^{2\alpha}}\right)^{\vee}(x)
\qe
(extended by 0 if $t\leq 0$ ) is the fundamental solution of the generalized heat equation, hence
\[
\big(\partial_t + (-\Delta_x)^\alpha \big) G_t = \delta
\]
in $D^\prime(\Rx)$.
Here $\wedge$ and $\vee$ stand for the Fourier transform and its inverse in $\Sd$, respectively and $\ast$ stands for the convolution in $S^\prime (\Rn)$.
\par
The paper is organized as follows. In Section 2 we fix notation, define function spaces and recall embeddings which are used to characterize mapping properties of $\Sm{t}$. The main results are contained in Section 3. Complements and comments are added in Section 4.

%===================================================================================================================================================
\section{Preliminaries}
\subsection{Function spaces}\label{fs}

Let $\Rn$ be the Euclidean $n\,$- space with $n\in\N$ where $\N$ indicates the collection of all natural numbers, $\N_0=\N\cup \{0\}$, 
\begin{eqnarray*}
	\Rxt=\{(x,t)\in\mathbb{R}^{n+1}:\;x\in\Rn,\;t>0\}
\end{eqnarray*} 
and $\overline{\Rxt}$ its closure. Put $\R=\R^1$. 
$C^u(\Rn)$ denotes the space of all $u\,$- times continuously differentiable functions in $\Rn$. $\Sw$ denotes the Schwartz space of all complex-valued infinitely differentiable rapidly decreasing functions on $\Rn$ and $\Sd$ its dual, the space of all tempered distributions. Furthermore, $\Lp{p}$ with $1\leq p\leq\infty$ is the Banach space of all $p\,$- integrable complex-valued functions with respect to the Lebesgue measure, normed by
\begin{eqnarray*}
	\|f|\Lp{p}\|=\left(\int_{\Rn}|f(x)|^p\text{d}x\right)^{1/p}
\end{eqnarray*}
with the usual modification if $p=\infty$.
If $\phi\in\Sw$ then
\begin{eqnarray}\label{Fou}
\phi^{\wedge}(\xi)=(\Fou\phi)(\xi)=(2\pi)^{-n/2}\int_{\Rn} e^{-ix\xi} \phi(x)\text{d}x,\quad \xi\in\Rn
\end{eqnarray}  
denotes the Fourier transform of $\phi$. $\Fou^{-1}\phi$ and $\phi^{\vee}$ stand for the inverse Fourier transform, given by the right hand side of \eqref{Fou} with $i$ in place of $-i$ and $x\xi$ stands for the scalar product in $\Rn$. $\Fou$ and $\Fou^{-1}$ are extended in the usual way to $\Sd$.
Let $\phi_0\in \Sw$ with 
\begin{eqnarray}\label{fi1}
\phi_0(x)=1 \text{ if } |x|\leq 1\; \text{ and }\;  \phi_0(x)=0 \text{ if } |x|\geq 3/2.
\end{eqnarray}  
We define the sequences
\begin{eqnarray}\label{fi2}
\phi_j(x)=\phi_0(2^{-j}x) -\phi_0(2^{-j+1}x),\quad x\in\Rn,\; j\in\N.
\end{eqnarray} 
Recall that $\Fou^{-1}[\phi_j\Fou f]$ and $\Fou^{-1}[\phi^j\Fou f]$ are entire analytic functions and make sense pointwise  for any $f\in\Sd$. We are interested in inhomogeneous Besov and Triebel-Lizorkin spaces $\A$, $A\in\{B,F\}$ with $s\in\R$, $1\leq p,q\leq\infty$ ($p<\infty$ for $F$-spaces). %($p<\infty$ for the $F$ -spaces).
The standard norms for these spaces are given by
%\begin{definition}\label{spaces}
\begin{eqnarray}\label{spaces}
\|f|\A\|=
\begin{cases}
\left(\sum_{j\in\N_0} 2^{jsq}\|\Fou^{-1}\phi_j\Fou f|\Lp{p}\|^q\right)^{1/q}, & \text{ if }A=B,\\
\|\left(\sum_{j\in\N_0} 2^{jsq}|\Fou^{-1}\phi_j\Fou f(\cdot)|^q\right)^{1/q}|\Lp{p}\|,& \text{ if }A=F
\end{cases}
\end{eqnarray} 
(with the usual modification if $q=\infty$). In particular the $\A$-spaces contain all tempered distributions $f\in\Sd$ such that the norm \eqref{spaces} is finite.
We recall that \eqref{spaces} is independent of the chosen resolution of unity (in the sense of equivalent norms). A detailed study of the inhomogeneous spaces including their history and properties can be found in \cite{Tri83}, \cite{Tri92} and \cite{Tri06}.\\ Finally, we recall that $B^s_{\infty,\infty}(\Rn)=\mathcal{C}^s(\Rn)$, $s\in\R$, in the sense of equivalent norms, where $\mathcal{C}^s(\Rn)$ denote the Hölder-Zygmund spaces. \\
As already sketched in the Introduction  we need  weighted $L_v$-spaces with respect to the Bochner integral. Let $X$ be a Banach space with $X\subset \Sd$, $0<T\leq\infty$, $b\in\R$ and $1\leq v \leq\infty$. Then $\Lv$ contains all $f:(0,T)\rightarrow X$ such that 
\begin{eqnarray}\label{Lv}
\|f|\Lv\|=\left(\int_0^T t^{bv}\|f(\cdot,t)|X\|^v\text{d} t\right)^{1/v} <\infty
\end{eqnarray}
if $1\leq v <\infty$ and
\eq\label{Lvinf}
\|f|\La{\infty}{b}{X}\|=\underset{0<t<T}{\sup}t^{b}\|f(\cdot,t)|X\| <\infty
\qe
if $v=\infty$.
We deal with Banach spaces $X=\A$, where $1\leq p,q\leq \infty$ ($p< \infty$ for $F$-spaces). \\
After extending functions satisfying either \eqref{Lv} or \eqref{Lvinf} from $\Rn\times(0,T)$ to $\mathbb{R}^{n+1}$ by zero one has always
\begin{eqnarray*}
	\La{v}{b}{\A}\subset S'(\mathbb{R}^{n+1})
\end{eqnarray*} 
if $s>0$ and $b<1-\frac{1}{v}$, cf. Remark \ref{distr}.
\par

%-------------------------------------------------------------------------------------------------------------------------------------------------
\subsection{Embeddings and multiplication properties}
In connection with the semi-linear generalized heat equation one has to deal with mapping properties of type $u\mapsto |u|^{r-1}u$ in the spaces $\A$. More precisely, we ask for conditions on the space $\A$ such that the product $|u|^{r-1}u$ with $u\in\A$ belongs again to $\A$. First recall that spaces $\A$ are algebras with respect to pointwise multiplication if $1\leq p,q\leq\infty$ and $s>\frac{n}{p}$. Then there exists some $c>0$ such that
\eqs
\|f_1f_2|\A\|\leq c\,\|f_1|\A\|\cdot \|f_2|\A\|
\qes %, i.e. $1<p<\infty$, $2\leq n\in\N$ and $s>\frac{n}{p}$. In particular $s>0$.\\
for all $f_1,f_2\in\A$, abbreviated as
\eqs
\A\cdot\A\hra\A.
\qes
The following Proposition summarizes some useful embeddings which we already adapted to our later needs.

\begin{prop}\label{LocLp}
	Let $1\leq p,q\leq\infty$ ($p<\infty$ for $F$-spaces), $s>0$ and $2\leq n\in\N$. Then
	\begin{itemize}
		\item [{\upshape (i)}]
		\eqs
		\A\hra\Lloc\quad\text{if, and only if,}\quad \A\hra\Lp{p}.
		\qes    	
		\item [{\upshape (ii)}] Let additionally be $r>1$ and $\frac{n}{p}<s<r$. Then it holds
		\eqs
		\||f|^r |\A\|\leq c\|f|\A\|^r.
		\qes
		\item [{\upshape (iii)}] Let additionally $f$ be real and in case of $F$-spaces $s\neq1$ if $p=1$. Then it holds
		\eqs
		\||f| |\A\|\leq c\|f|\A\|\quad\text{if, and only if,}\quad 0<s<1+\frac{1}{p}.
		\qes
	\end{itemize}
\end{prop}	
\begin{proof}
	Part (i) follows from \cite[Theorem 3.3.1, Theorem 3.3.2 and Corollary 3.3.1]{ST}. Part (ii) is a consequence of  \cite[Theorem 1, Chapter 5.4.3]{RuS} and the embedding $\A\hra\Lp{\infty}$. Part (iii) is also known as trunction property and follows from \cite[Theorem 25.8]{Tri01}.
\end{proof}
\begin{rem}\label{mü}
	Part (ii) also holds if one replaces $|f|^r$ by $|f|^{r-1}f$. This follows from \cite[Remark 4, Chapter 5.4.3]{RuS}
\end{rem}
\begin{rem}
	In case of Besov spaces assertion (ii) in Proposition \ref{LocLp} even holds for $\frac{n}{p}<s<r+\frac{1}{p}$. This was shown in \cite{K}.
\end{rem}
\section{Semi-linear hyperdissipative heat equations}
Due to the biological background we assume that initial data $u_0$ in \eqref{ghsup} are real valued. From the structure of the problem it follows that then solutions $u$ are also real. Therefore, we deal from the very beginning with the real parts of $\A$, including related fixed point assertions.
First we recall
recall \cite[Theorem 3.5]{BaS19}. The associated estimates are essential for the proof of Theorem \ref{fix2}.
 
\begin{prop}\label{esthom}
	Let $1\leq p,q\leq\infty$, $\alpha\in\N$,  $s\in\R$, $\alpha\in\N$ and $d\geq0$. Then there is a constant $c>0$ such that 
	\eq\label{firstest}
	\|W_{t}^{\alpha}\omega|A^{s+d}_{p,q}(\Rn)\|\leq c\, t^{-d/2\alpha} \|\omega|\A\|
	\qe
	for all $t$ with $0<t\leq 1$ and $\omega\in\A$.
\end{prop}
As already indicated we consider solutions of \eqref{ghsup} which are fixed points of the related integral operator \eqref{T0}. Such a solution is called mild. For more detailed explanations we refer for instance to \cite{Can} and \cite{Lem02}. We rely on \eqref{T0} in the context of the spaces $\A$. Moreover, we ask in which spaces $u(x,t)$ converges to the initial data $u_0$ if $t$ tends to zero, hence in which spaces the solution is strong. Before we deal with the fixed point problem, we need some a priori estimates of $\T$ where we use the above Proposition.

\begin{prop}\label{estinhom}
	Let $2\leq n\in\N$, $\alpha\in\N$, $1\leq p,q\leq\infty$ ($p<\infty$ for $F$-spaces), $\frac{n}{p}<s$ and $\max\{1,s\}<r$. Let $T>0$ and $a$, $v$ such that
	\eq\label{Bed}
	\frac{1}{2}\leq v\leq\infty\quad\text{ and }\quad -\infty<a<2- \frac{1}{v}.
	\qe
	Let 
	\eq
	u_0\in \As{s_0}\mbox{ with } s_0\leq s\mbox{ and } u\in\La{2r v}{\tfrac{a}{2r}}{\A}.
	\qe
	Then there is a constant $c>0$, independent of $u_0$ and $u$, such that 
	\eq\label{Tuo}
	\|\T u(\cdot, t)|\A\|&\leq& c\,t^{-\frac{s-s_{0}}{2\alpha}} \|u_0|\As{s_0}\|\\
	&+&c\,t^{1-\frac{a}{2}-\frac{1}{2v}}
	\|u|\La{2r v}{\tfrac{a}{2r}}{\A}\|^r.\notag
	\qe
	for all $t$ with $0<t<T$ (with $\frac{1}{v}=0$ and the modification \eqref{Lvinf} if $v=\infty$).
\end{prop}
\begin{proof}
	Concerning the first summand we apply Proposition \ref{esthom} with $d=s-s_0$. For the second summand we obtain with $d=0$ and Remark \ref{mü}
	\begin{align*}
	& \int_0^t\|\Sm{t-\tau}|u|^{r-1}u(\cdot,\tau)|\A\|\,\text{d}\tau\
	\leq c \int_0^t \||u|^{r-1}u(\cdot,\tau)|\A\|\,\text{d}\tau\\
	&\leq c \int_0^t\tau^{-\frac{a}{2}}\tau^{\frac{a}{2}}\,\|u(\cdot,\tau)|\A\|^r\,\text{d}\tau\\
	&\leq c
	\left(\int_0^t\tau^{-\frac{a}{2}(2 v)'}\text{d}\tau\right)^{\frac{1}{(2 v)'}}\left(\int_0^t \tau^{av}\|u(\cdot,\tau)|\A\|^{2r v}\,\text{d}\tau\right)^{\frac{1}{2v}}\\
	&\leq c\,t^{1-\frac{a}{2}-\frac{1}{2v}}\,\|u|\La{2rv}{\tfrac{a}{2r}}{\A}\|^r.
	\end{align*}
\end{proof}

\begin{theorem}\label{fix2}
	Let $2\leq n\in\N$, $\alpha\in\N$, $1\leq p,q\leq\infty$  ($p<\infty$ for $F$-spaces), $s>\frac{n}{p}$ and $r\geq2$. 
	\bli
	\item[{\upshape\bfseries (i)}]Let $u_0\in\As{s_0}$ be real valued initial data with
	\begin{align}
	\frac{n}{p}-\frac{2\alpha}{r-1}<s_0\leq s<r-1.\label{e_1}
	\end{align}
	Then there exists a number $T>0$ such that \eqref{ghsup} has a unique real mild solution 
	\begin{eqnarray*}
		u\in\La{2rv}{\tfrac{a}{2r}}{\A} \cap  \Cinf
	\end{eqnarray*}
whereas $a$ and $v$ fulfill the conditions
	\begin{align}
	&\frac{1}{2}< v\leq\infty \label{e_2}\\
	&\frac{r(s-s_0)}{\alpha}<a+\frac{1}{v}<2\label{e_3}. 
	\end{align}

	\item[{\upshape\bfseries (ii)}]
	If, in addition, $\max\{p,q\}<\infty$ then the above solution is strong, that means $u\in C([0,T),\As{s_0})$.
	\eli
\end{theorem}
\begin{proof}
	\underline{\textbf{Step 1.}}
	Due to Proposition \ref{estinhom} we have
	\begin{align*}
	\|T_{u_0}u(\cdot,t)|&\A)\|\leq c\, 
	t^{-\frac{s-s_{0}}{2\alpha}} \|u_0|\As{s_0}\|\\
	&+c\,t^{1-\frac{a}{2}-\frac{1}{2v}}
	\|u|\La{2rv}{\tfrac{a}{2r}}{\A}\|^r.
	%\label{s10}
	\end{align*}
	for $0<t<T$.
	We multiply both sides with $t^{\frac{a}{2r}}$. Raising to the power of $2rv$ and integrating over $(0,T)$ yields
	\begin{align}\label{s11}
	&\int_0^T t^{av}\|T_{u_0}u(\cdot,t)|A^s_{p,q}(\Rn)\|^{2r v}\text{d}t\\
	&\leq  c\,T^{\delta}\,\|u_0| \As{s_0}\|^{2r v}\nonumber
	+T^{\varkappa}\,\|u|\La{2r v}{\tfrac{a}{2r}}{\A}\|^{2vr^2}
	\end{align}	
	with
	\begin{align}\label{s12a}
	\delta=av-\frac{r(s-s_0)}{\alpha}v+1>0,~~\mbox{since }~ \frac{r(s-s_0)}{\alpha}<a+\frac{1}{v}
	\end{align}
	and 
	\begin{align}\label{s12b}
	\varkappa=av+2rv-arv-r+1>0,~~\mbox{since } ~a+\frac{1}{v}<2<\frac{2r}{r-1}.
	\end{align}
	Moreover, condition \eqref{e_1} ensures that \eqref{e_3} makes sense.
	Thus, $T_{u_0}$ maps the unit ball $U_T$ in $\La{2rv}{\tfrac{a}{2r}}{\A}$ into itself if $T$ is sufficiently small.\\
	\\
	As for the contraction we consider a function $f\in C^1(\R)$, $-\infty<x_1<x_2<\infty$ and $y=x_1+t(x_2-x_1)$, $t\in[0,1]$. Then we can write
	\begin{align*}
	f(x_2)&=f(x_1)+\int\limits_{x_1}^{x_2} f'(y)\text{d}y\\
	&=f(x_1) +(x_2-x_1)\int\limits_0^1 f'(tx_2+(1-t)x_1)\,\text{d} t.
	\end{align*}
	This applies to the function $f(x)=|x|^{r-1}x$ with $x\in\R$ and $r>1$.	Then we have
	\begin{align*}
	|x_2|^{r-1}x_2=|x_1|^{r-1}x_1+r(x_2-x_1)\int\limits_0^1 \big|tx_2+(1-t)x_1\big|^{r-1}\text{d} t.
	\end{align*}
	With $x_2=u(x)\geq v(x)=x_1$ one obtains the identity
	\begin{align}\label{s13}
	|u(x)|^{r-1}u(x)&-|v(x)|^{r-1}v(x)\\
	&= r\big(u(x)-v(x)\big)\int\limits_0^1 |tu(x)+(1-t)v(x)|^{r-1}\,\text{d} t\nonumber
	\end{align}
	for $x\in\Rn$ and $r>1$. If $v(x)>u(x)$ real then $u$ and $v$ change their values but the outcome is the same.\\
	\\
	Let now belong $u,v$ to $U_T$, $T>0$. Using the fact that $\A$ are multiplication algebras and Banach spaces and Proposition \ref{LocLp} part (ii) with $r-1$ instead of $r$, $2<r$ and part (iii) if $r=2$ then we obtain
	\begin{align}\label{s14}
	&\|\left(|u|^{r-1}u\right)(\cdot,\tau)-\left(|v|^{r-1}v\right)(\cdot,\tau)|\Asp\|\notag\\
	&\leq c \|u(\cdot,\tau)-v(\cdot,\tau)|\Asp\|\cdot \int\limits_0^1\left\||tu(\cdot,\tau)+(1-t)v(\cdot,\tau)|^{r-1}\big|\Asp\right\|\text{d} t\nonumber\\
	&\leq c \|u(\cdot,\tau)-v(\cdot,\tau)|\Asp\|\cdot \int\limits_0^1 \left(t\|u(\cdot,\tau)|\Asp\|+(1-t)
	\|v(\cdot,\tau)|\Asp\|\right)^{r-1}\text{d} t\nonumber\\
	&\leq c_r \|u(\cdot,\tau)-v(\cdot,\tau)|\Asp\|\cdot\Big(\|u(\cdot,\tau)|\Asp\|^{r-1}+\|v(\cdot,\tau)|\Asp\|^{r-1}\Big)
	\end{align}
	where we abbreviated $\A$ by $\Asp$ for better readability. It is sufficient to calculate the estimate for the first summand. The second one can be treated similarly.
	Applying twice Hölder's inequality and using the estimate in \eqref{s11} we obtain
	\begin{align}
	&\|\T u - \T v|\La{2rv}{\tfrac{a}{2r}}{\A}\|^{2rv}\\
	&\leq c\int\limits_0^T t^{av}\left(\int\limits_0^t  \left\|\left(|u|^{r-1}u\right)(\cdot,\tau)-\left(|v|^{r-1}v\right)(\cdot,\tau)|\A\right\| \text{d}\tau\right)^{2rv}\text{d} t\nonumber\\
	%& \leq c\int\limits_0^T t^{av}\left(\int\limits_0^t \left\|u(\cdot,\tau)-v(\cdot,\tau)|\Asp\|\cdot\|u(\cdot,\tau)|\Asp\right\|^{r-1}\text{d}\tau\right)^{2rv}\text{d} t\nonumber\\
	&\leq cT^{\varkappa}\|(u-v)|\La{2rv}{\tfrac{a}{2r}}{\A}\|
	\end{align}
	with some $\varkappa>0$. Hence $\T: U_T\to U_T\subset \La{2rv}{\tfrac{a}{2r}}{\A}$ is a contraction if $T>0$ is sufficiently small. Application of Banach's contraction principle yields that $\T$ has a unique fixed point in the unit ball of $\La{2rv}{\tfrac{a}{2r}}{\A}$ which is a solution of \eqref{ghsup}.\\
	\underline{\textbf{Step 2.}} 
	Before we are able to extend the above result from the unit ball $U_T$ to the whole space $\La{2rv}{\tfrac{a}{2r}}{\A}$ we need some technical preparation. First we show that under the above conditions the fixed point $\T u(x,t)=u(x,t)$ is a $C^\infty$ -function in the open strip $\R^n\times (0,T)$. Similarly as in the proof of Proposition \ref{estinhom} we obtain with $s+\eta$ instead of $s$ and $d=s+\eta-s_0$ and $d=\eta$ in the estimate the first and second summand of $\eqref{T0}$, respectively
	\begin{align}
	\|T_{u_0}u(\cdot,t)&|\As{s+\eta}\|\leq c\,t^{-\frac{s+\eta-s_0}{2\alpha}}\|u_0|\As{s_0}\|\notag\\
	&+c\,\int_0^t (t-\tau)^{-\frac{\eta}{2\alpha}}\||u|^{r-1}u(\cdot,\tau)|\A\|\text{d}\tau\notag\\
	&\leq c\,t^{-\frac{s+\eta-s_0}{2\alpha}}\|u_0|\As{s_0}\|\label{1}\\
	&+c\, t^{1-\frac{1}{2v}-\frac{\eta}{2\alpha}-\frac{a}{2}}\|u|\La{2rv}{\tfrac{a}{2r}}{\A}\|^r.\label{2}
	\end{align}
	Here $\eta>0$ must be chosen such that $\max(1,s+\eta)<r$ and $-\infty <a+\frac{\eta}{2}<2-\frac{1}{v}$.\\ We iterate this argument with smoother initial data $u_\varepsilon(x)=u(x,\varepsilon)$ where $u$ is the solution obtained in Step 1 and $\varepsilon>0$ is chosen sufficiently small. Then one has $u_\varepsilon(x)\in \As{s+\eta}\hra\As{s_0+\eta}$. Hence
	\begin{align}
	\|T_{u_\varepsilon}u(\cdot,t)&|\As{s+2\eta}\|\leq c\,t^{-\frac{s+\eta-s_0}{2\alpha}}\|u_\varepsilon|\As{s_0+\eta}\|\notag\\
	&+c\,\int_0^t (t-\tau)^{-\frac{\eta}{2\alpha}}\||u|^{r-1}u(\cdot,\tau)|\As{s+\eta}\|\text{d}\tau\notag\\
	&\leq c\,t^{-\frac{s+\eta-s_0}{2\alpha}}\|u_\varepsilon|\As{s_0+\eta}\|\label{3}\\
	&+c\, t^{1-\frac{1}{2v}-\frac{\eta}{2\alpha}-\frac{a}{2}}\|u|\La{2rv}{\tfrac{a}{2r}}{\As{s+\eta}}\|^r.\label{4}
	\end{align}
	In particular we see that the factors in \eqref{1} and \eqref{2} do not change. Thus, after the $k$-th iteration step we arrive at
	\begin{align}
	\|T_{u_\varepsilon}u(\cdot,t)&|\As{s+k\eta}\|\leq c\,t^{-\frac{s+\eta-s_0}{2\alpha} }\|u_\varepsilon|\As{s_0+(k-1)\eta}\|\notag\\
	&+c\, t^{1-\frac{1}{2v}-\frac{\eta}{2\alpha}-\frac{a}{2}}\|u|\La{2rv}{\tfrac{a}{2r}}{\As{s+(k-1)\eta}}\|
	\end{align}
	The embedding $\AsP{s+k\eta}(\Rn)\hra B^{s+k\eta-\frac{n}{p}}_{\infty,\infty}(\Rn)=\mathcal{C}^{s+k\eta-\frac{n}{p}}(\Rn)$ for any $k\in\N$ yields that the solution $u$ is infinite differentiable with respect to the space variable. Note that the above arguments do not affect the conditions on $a$ and $v$. \\We incorporate the time variable. From Step 1 we have $u\in\La{\infty}{\frac{a}{2r}}{\A}$ and hence,
	\eqs
	\underset{t\in(t_0,t_1)}{\sup}t^{\frac{a}{2r}}\|u(\cdot,t)|\A\|<\infty,\quad 0<t_0<t_1<T.
	\qes
	This leads to 
	\begin{align*}
	\underset{t\in (t_0,t_1)}{\sup}\|u(\cdot,t)|\A\|&=\underset{t\in (t_0,t_1)}{\sup}t^{-\frac{a}{2r}}t^{\frac{a}{2r}}\|u(\cdot,t)|\A\|\\
	&\leq
	t_0^{-\frac{a}{2r}}\underset{t\in (t_0,t_1)}{\sup} t^{\frac{a}{2r}}\|u(\cdot,t)|\A\|
	<\infty~
	\end{align*}
	and in combination with the above considerations to
	\begin{align}\label{dx}
	\underset{x\in\Omega,\,t\in (t_0,t_1)}{\sup}|\text{D}_x^{\beta}u(x,t)|\leq c_{\beta}\quad \text{for all } \beta\in\N_0^n.
	\end{align}
	where $\Omega\subset\Rn$ is a bounded domain. Moreover, $u$ solves
	\eq\label{dt}
	\partial_t u+(-\Delta_x)^\alpha u=|u|^{r-1}u
	\qe
	in the distributional sense, see Chapter 4, Remark \ref{distr}. Because of \eqref{dx} the RHS of \eqref{dt} is uniformly bounded in $\Omega\times(t_0,t_1)$, $0<t_0<t_1<T$ and thus $\partial_t u$. Iterating this argument and using \eqref{dx} we finally have
	\begin{align*}
	\underset{x\in\Omega,\,t\in (t_0,t_1)}{\sup}\Big|\frac{\partial^k}{\partial t^k} \text{D}_x^{\beta}u(x,t)\Big|\leq c_{\beta},\quad \text{for all } \beta\in\N_0^n,
	\end{align*}
	hence, 
	$
	\frac{\partial^k}{\partial t^k} \text{D}_x^{\beta}u\in L_{\infty}(\Omega\times (t_0,t_1))
	$. Application of Sobolevs embedding theorem leads to the desired result for any $k\in\N$, $\beta\in\N_0$, i.e. $u\in C^\infty(\Rn\times(0,T))$.\\
	\underline{\textbf{Step 3.}}
	Let temporarily denote $X_T:=\La{2rv}{\tfrac{a}{2r}}{\A}$. Assume moreover, that $u\in U_T$ is the solution obtained in Step 1 and let $v\in X_T$ be a second solution of \eqref{ghsup}. Then by the same arguments as in Proposition \ref{estinhom} and \eqref{s14} we have
	\eqs
	\|\T u-\T v|X_{T_0}\|\leq c\, T_0^{\varkappa}\|u-v| X_{T_0}\|\cdot(1+\|v|X_T\|^{r-1}), 
	\qes
	$0<T_0\leq T.$ We choose $T_0$ such that $c\, T_0^{\varkappa}\cdot(1+\|v|X_T\|^{r-1})<1$ and thus it holds $u(\cdot, t)=v(\cdot,t)$ for almost all $0<t\leq T_0$. Now one can apply the arguments in \cite[Theorem 5.24 Step 2, Remark 5.25]{Tri13} to show that the solution $u$ is not only unique in $U_T\subset X_T$ but in the whole space.
	\\
	\underline{\textbf{Step 4.}} We show that the solution $u(x,t)$ is also strong if $\max\{p,q\}<\infty$. Therefore, we need to show that the solution converges to the initial data $u_0$, if $t\to 0$ in the \ $\As{s_0}$ - norm. We obtain
	\begin{align}\label{cont}
	&\|u(\cdot,t)-u_0|\As{s_0}\|\notag\\
	&\leq \|\Sm{t}u_0-u_0|\As{s_0}\|
	+c\int_0^t \|\Sm{t-\tau}\left(|u|^{r-1}u\right)(\cdot,\tau)|\As{s_0}\|\,\text{d}\tau.
	\end{align}
	The second summand in \eqref{cont} we estimate from above in the $\A$ - norm, since $s\geq s_0$ and use Proposition \ref{estinhom}. This yields 
	\begin{align}\label{est-1}
	\int_0^t \|\Sm{t-\tau}\left(|u|^{r-1}u\right)(\cdot,\tau)&|\As{s_0}\|\,\text{d}\tau\\ & \leq c\,t^{1-\frac{a}{2}-\frac{1}{2v}}
	\|u|\La{2r v}{\tfrac{a}{2r}}{\A}\|^r,\notag
	\end{align}
	if $v<\infty$. The right hand side tends to zero for $t$ tending to zero since $a<2-\frac{1}{v}$. If $v=\infty$ we use the modification \eqref{Lvinf}. Then $a<2$ is required to ensure convergence. Concerning the first summand one can proceed similarly to \cite [Theorem 3.8, Step 5]{BaS17}.
\end{proof}
\begin{rem}
	Note that if $r\in\N$ is an odd number then the restriction on the smoothness $s$ is not necessary. Moreover, one can ask for solutions belonging to function spaces beyond multiplication algebras.
\end{rem}
In addition to the results of the previous part one may ask for stability of solutions. That means small perturbations of the initial data cause small deviations of the solution. Solutions $u_1$ and $u_2$ of \eqref{ghsup} according to Theorem \ref{fix2} with respect to initial data $u_0^1$ and $u_0^2$, respectively, 
are called locally stable if for any $\varepsilon>0$ there exists a $\delta>0$ and a time $T>0$ such that for all $0<t<T$ holds
\eq\label{eps}
\|u_1(\cdot,t)-u_2(\cdot,t)|\As{s_0}\|\leq \varepsilon
\qe
if
\eq\label{del}
\|u^1_0-u^2_0|\As{s_0}\|\leq \delta~.
\qe
A problem is called well-posed if there are unique mild solutions which are additionally strong and stable in the above sense.
It is sufficient to show the stability for
$u\in\La{\infty}{\frac{a}{2r}}{\A}$. Recall that by construction of the solution as a fixed point of $\T$ we have $\|u|\La{\infty}{\frac{a}{2r}}{\A}\|\leq 1$.
\par

\begin{theorem}\label{stablesol}
	Let $u_i\in L_\infty ((0,T_i),\,\frac{a}{2r},\,\A )$ ($i=1,2$) be solutions of \eqref{ghsup} obtained in Theorem \ref{fix2} 
	with initial data $u^i_0\in\As{s_0}$ in the corresponding time interval $(0,T_i)$. 
	Then under the conditions of Theorem \ref{fix2} Part (i) the solutions are locally stable. 
	If additionally $\max \{p,q\}<\infty$ then the Cauchy problem \eqref{ghsup} is well-posed.
\end{theorem}
\begin{proof}
	Let $u_1$, $u_2$ be two solutions of \eqref{ghsup}, with corresponding initial data $u_0^1$, $u_0^2$. Then it holds
	\begin{align}
	& \|u_1(\cdot,t)-u_2(\cdot,t)|\As{s_0}\| 
	\leq  \|\Sm{t}(u^1_0-u^2_0)|\As{s_0}\|\notag\\
	&+\int_0^t\|\Sm{t-\tau}(|u_1|^{r-1}u_1-|u_2|^{r-1}u_2)(\cdot,\tau)|\As{s_0}\|\text{d}\tau\notag\nonumber\\
	%&\text{Since\ \  $s>s_0 \implies \A\subset \As{s_0}$; }\nonumber\\
	&\leq c_1\|(u^1_0-u^2_0)(\cdot)|\As{s_0}\|\notag\\
	&+c_2\int_0^t\|\Sm{t-\tau}(|u_1|^{r-1}u_1-|u_2|^{r-1}u_2)(\cdot,\tau)|\A\|\text{d}\tau\notag\\
	&\leq c_1\delta
	+c_2\int_0^t\|(|u_1|^{r-1}u_1-|u_2|^{r-1}u_2)(\cdot,\tau)|\A\|\text{d}\tau\notag,
	\end{align}
	where we used elementary embedding and Proposition \ref{esthom} with $d=0$. Hence, using the same idea as for the proof of the contraction now with $u_1(\cdot,t)\geq u_2(\cdot,t)$ we obtain
	\begin{align}
	&\|u_1(\cdot,t)-u_2(\cdot,t)|\As{s_0}\| \leq c_1\delta\notag\\
	&+c_2\sum\limits_{i=1}^2\int_0^t\|(u_1-u_2)(\cdot,\tau)|\A\|\|u_i(\cdot,\tau)|\A\|^{r-1}\text{d}\tau\label{3.17}%\notag
	\end{align}
	By means of Minkowski's inequality we split the norm in the second summand in the right part of the inequality \eqref{3.17} and obtain for $i=1$
	\begin{align}
	&\int_0^t\|(u_1-u_2)(\cdot,\tau)|\A\|\|u_1(\cdot,\tau)|\A\|^{r-1}\text{d}\tau\notag\\
	&\leq \int_0^t\tau^{-\frac{a}{2}}\,\text{d}\tau\left(\|u_1|X_\infty^1\|^r+\|u_1|X_\infty^1\|^{r-1}\|u_2|X_\infty^2\|\right)\notag
	\end{align}
	where we abbreviated $X^i_\infty:=L_\infty ((0,T_i),\,\frac{a}{2r},\,\A )$. Finally, 
	\begin{align}
	\|u_1(\cdot,t)-u_2(\cdot,t)|\As{s_0}\| \leq c_1\delta
	+ct^{-\frac{a}{2}+1}
	\end{align}
	since $\|u_i|X^i_\infty\|\leq 1,$ $i=1,2$.
	Due to condition \eqref{Bed} the term $t^{-\frac{a}{2}+1}$ tends to $0$ for $t\to 0$. Hence,
	\begin{align}
	\|u_1(\cdot,t)-u_2(\cdot,t)|\As{s_0}\|\leq \varepsilon
	\end{align}
	for all $0<t<T$ if $T\leq \min\{T_1,T_2\}$ is small enough.
	Thus, the solution is locally stable and strong if additionally $\max\{p,q\}<\infty$ and the Cauchy problem \eqref{ghsup} well-posed.
\end{proof}

\section{Complements and comments}
We add some comments on uniqueness and the spaces we are dealing with.  
\begin{rem}\label{distr}
Theorem \ref{fix2} does not answer the question of uniqueness of solutions of \eqref{ghsup} in the framework of $D'(\Rx)$. Nevertheless, the  uniqueness result proved there  mainly relies on the fact that the fixed point of the operator $\T$ is also a solution of \eqref{ghsup} in the distributional sense, see for instance \cite[Sections 3.3.4-3.3.6]{Tri92}. Note at first that under the conditions in Theorem \ref{fix2} one has always
\eqs
\La{2rv}{\tfrac{a}{2r}}{\A}\subset D'(\Rx)
\qes
after extending functions belonging to this space from $\Rn\times (0,T)$ to $\Rx$ by zero. This can be seen using a duality argument. Let $\varphi(x,t)\in D(\Rx)$ and $u(x,t)\in \La{2rv}{\frac{a}{2r}}{\A}$. Using 
\eqs
\Big(B^{-s+\frac{n}{p}}_{1,1}\Big)(\Rn)'=B^{s-\frac{n}{p}}_{\infty,\infty}(\Rn)=\mathcal{C}^{s-\frac{n}{p}}(\Rn)
\qes
we obtain
\begin{align}
|\langle u(\cdot,t),\varphi(\cdot,t)\rangle|&\leq c\,\|u(\cdot,t)|\mathcal{C}^{s-\frac{n}{p}}(\Rn)\|\|\varphi(\cdot,t)|B^{-s+\frac{n}{p}}_{1,1}(\Rn)\|\\
&\leq c'\,\|u(\cdot,t)|\A\|,
\end{align}
where $c'$ is independent of $t$. Integration then yields
\begin{align}
|\langle u,\varphi\rangle|&\leq c'\int\limits_0^T\|u(\cdot,t)|\A\|\text{d}t\\
& \leq c'\Bigg(\int\limits_0^T t^{-\frac{a}{2r}(2rv)'}\text{d}t\Bigg)^{1/(2rv)'} \|u|\La{2rv}{\tfrac{a}{2r}}{\A}\|.
\end{align}
If $\frac{a}{2r}(2rv)'<1$ then all spaces we are dealing with makes sense in $D'(\Rx)$ and even in $S'(\Rx)$ (similar argument with $\varphi(x,t)\in S(\Rx)$).  Moreover, one can show that under the above conditions $|u|^{r-1}u$ and the RHS of \eqref{T0} is locally integrable on $\Rx$ (again after extending functions by zero). Here one can use Proposition \ref{estinhom} and follow the ideas e.g. in \cite[Lemma 3.5]{BaS19}.
 Since $\K{t}$ as defined in \eqref{K} is a fundamental solution of \eqref{ghsup} and $\Sm{t}u_0$ solves the homogeneous equation in the classical sense, see again \cite[Lemma 3.5]{BaS19}, we can deduce with \cite[Theorem 3.2.4/2]{T92} that any mild solution of \eqref{ghsup} is also a solution in the distributional sense. 
\end{rem}
We add a remark on (sub/super)-critical spaces and the lower bound in \eqref{e_1}.
\begin{rem}\label{rem-3.10}
	As far as initial data are concerned one can ask for critical, sub- and supercritical spaces. We follow the terminology used e.g. in \cite[Subsection 6.2.5/3, p.211]{Tri13}, \cite[Chapter 2]{Tri17} or \cite[Subsection 3.2]{BaS19} and the references given therein. The approach is based on scaling properties of the solution of \eqref{ghsup}. Let $\lambda>0$ and $u(x,t)$ be a solution of \eqref{ghsup} in $\Rn\times(0,T)$ with initial data $u_0$. Then we easily observe that
	\begin{align}
	u_{\lambda}(x,t)&=\lambda^{\frac{2\alpha}{r-1}}u(\lambda x,\lambda^{2\alpha} t)
	\end{align}
	solves \eqref{ghsup} on $\Rn\times (0,\lambda^{-2\alpha}\,T)$ associated with the initial data
	\begin{align}
	u_{\lambda}(x,0)&=\lambda^{-\frac{2\alpha}{r-1}}u_0(\lambda^{-1}x).
	\end{align}
For the purpose of this consideration we assume that $u_0$ belongs to an homogeneous space	$\dot{A}^{s_0}_{p,q}(\Rn )$ (see \cite[Chapter 5]{Tri92} for related norms). Then it holds
\eqs
\|f(\lambda\cdot)|\dot{A}^{s_0}_{p,q}(\Rn )\|=\lambda^{s-\frac{n}{p}}\|f|\dot{A}^{s_0}_{p,q}(\Rn )\|, f\in \dot{A}^{s_0}_{p,q}(\Rn ),\;\lambda>0.
\qes
Let us suppose that there exists some $\delta >0$ such that for all $u_0\in\dot{A}^{s_0}_{p,q}(\Rn )$
with 
\[
\|u_0|\dot{A}^{s_0}_{p,q}(\Rn)\| \leq\delta
\]
there exists a solution $u(x,t)$ of \eqref{ghsup}. Using the above scaling argument and the homogeneity property of 
$\dot{A}^{s_0}_{p,q}(\Rn )$ we find that for all initial data $u_0\in\dot{A}^{s_0}_{p,q}(\Rn )$ with
\[
\|u_0|\dot{A}^{s_0}_{p,q}(\Rn)\| <\delta \lambda ^{s_0- \frac{n}{p}+\frac{2\alpha}{r-1}}
\]
there exists also a solution of \eqref{ghsup} on $\Rn\times (0,\lambda ^{-2\alpha}T)$.\\
If $s_0=\frac{n}{p} -\frac{2\alpha}{r-1}$ this leads to the existence of a  global solution if one considers $\lambda \to 0$.
This case is usually called critical. On the other hand we can expect the existence of local (small $T$) solutions
for arbitrarily large initial values in the supercritical case $s_0>\frac{n}{p} -\frac{2\alpha}{r-1}$. Whereas in the subcritical case 
$s_0<\frac{n}{p} -\frac{2\alpha}{r-1}$  the existence of global solutions cannot be expected.\\
The above scaling heuristics helps to get an imagination which (super)critical spaces are naturally adapted to problem \eqref{ghsup}. But it does not yield sharp boundaries to other spaces. In particular, our considerations rely on inhomogeneous spaces $\A$. In Theorem \ref{fix2} it comes out that our solution spaces belong to the supercritical area. This follows from condition \eqref{e_3} and is indicated in the lower bound for $s_0$ in \eqref{e_1}.
\end{rem}	

\end{document}